\documentclass{amsart}
\usepackage[utf8]{inputenc}
\usepackage{amsthm}
\usepackage{amsfonts}
\usepackage{amsmath}
\usepackage{mathtools}
\usepackage{amssymb}
\usepackage{url}
\usepackage{fullpage}
\usepackage{comment}
\usepackage{commath}
\usepackage{graphicx}
\usepackage{subcaption}
\usepackage{standalone}
\usepackage{tikz}
\usepackage{pgfplots}
\usepackage[parfill]{parskip}
\usepackage{longtable}
\usepackage[T1]{fontenc}
\usepackage{currvita}
\pgfplotsset{compat=1.17}
\usepackage{thm-restate}
\usepackage{listings}
\usepackage{array}
\usepackage[OT2,T1]{fontenc}
\usepackage[backend=biber, style=alphabetic]{biblatex}
\usepackage{xcolor}
\usepackage[colorlinks = true,
            linkcolor = blue,
            urlcolor  = blue,
            citecolor = blue,
            anchorcolor = blue]{hyperref}
\newcommand{\MYhref}[3][blue]{\href{#2}{\color{#1}{#3}}}%
\DeclareSymbolFont{cyrletters}{OT2}{wncyr}{m}{n}
\DeclareMathSymbol{\Sha}{\mathalpha}{cyrletters}{"58}

\makeatletter
\def\subsubsection{\@startsection{subsubsection}{3}%
  \z@{.5\linespacing\@plus.7\linespacing}{.1\linespacing}%
  {\normalfont\itshape}}
\makeatother

\title{Sphere Packing Densities of Certain Families of Elliptic Curves}
\author{Arjun Nigam}
\date{\today}

\theoremstyle{plain}
\newtheorem{thm}{Theorem}[section]
\newtheorem*{thm*}{Theorem}

\newtheorem{prop}[thm]{Proposition}

\theoremstyle{definition}
\newtheorem{defin}[thm]{Definition}

\theoremstyle{remark}
\newtheorem*{rem}{Remark}

\numberwithin{equation}{section}

\usepackage{indentfirst}
\begin{document}
\maketitle

\begin{abstract}
    In this paper, we examine the Mordell-Weil lattices of two families of elliptic curves over fields of characteristic $p>0$. We compute explicit lower bounds on the densest sphere packings of Mordell-Weil lattices using geometric information about the curve such as its Néron–Tate height and rational points on the curve. We see how these methods may be generalized to other ``nice" families of elliptic curves and how encode interesting properties about the elliptic curves themselves.
\end{abstract}
\section{Introduction}
A sphere packing in $\mathbb{R}^n$ is an arrangement of non-overlapping spheres of the same radius. The densest sphere packing in $n$ dimensions is a sphere packing in $\mathbb{R}^n$ such that the fraction of the volume occupied by the spheres is maximal. Lagrange proved in $1773$ that the densest sphere packing in two dimensions is the hexagonal packing. More than 200 years later, Thomas Hales proved that the cubic close packing and hexagonal close packing arrangements achieve the densest sphere packing in three dimensions. Later, Maryna Viazovska resolved the problem of finding the densest sphere packing for dimension $8$ in \cite{Via}. Together with Henry Cohn, Abhinav Kumar, Stephen Miller, and Danylo Radchenko, she later applied her method in \cite{CKMRV} to resolve the problem in dimension $24$ as well.

In this paper, we find lower bounds on the sphere packing densities associated to Mordell-Weil lattices of two families of elliptic curves. After introducing the relevant definitions and general results of sphere packings, we relate the density of the Mordell-Weil lattice of an elliptic curve to its Néron–Tate height function and regulator in section 3. 

In section 4, we focus on a specific family of elliptic curves and obtain lower bounds on their densities in two different ways. These bounds are obtained using subtle properties of the given family of elliptic curves. We compare our results to the best known sphere densities in \cite{Coh}. We also see how the information about the two lower bounds can be combined to obtain results about the original elliptic curve, namely a lower bound on the order of its Tate–Shafarevich group.  

Next, we focus on a different family of elliptic curves, one where we have more limited information. While we are unable to obtain a lower bounds on the densities of the Mordell-Weil lattices, we are able to find lower bounds on the densities of certain sublattices. We see that while this method requires less information about the elliptic curves than the method employed in the previous section, the resulting lower bounds on densities turn out to be much smaller than the best known sphere packing densities in the corresponding dimensions. 

In section 6, we explore how the ideas of the previous two sections may be applied to other families of elliptic curves. It is often the case that subtle properties of specific families of elliptic curves allow us to deduce results that do not hold for elliptic curves in general. However, there are many examples in the literature where the same result can be proven for different families of elliptic curves using different arguments. We state a few examples relevant to this paper and describe how the methods of this paper may be modified to study the sphere packing densities associated to other families of elliptic curves.
\section{Preliminaries}
We now make some basic definitions and introduce notation that we will use throughout the paper. 
\begin{defin}
    Let $\Lambda$ be a free $\mathbb{Z}$-module in $\mathbb{R}^n$ of maximal rank. A sphere packing on $\Lambda$ is a non-overlapping arrangement of spheres of the same radius in $\mathbb{R}^n$ such that each sphere is centered at a point of $\Lambda$. The \textbf{densest sphere packing} of $\Lambda$ is a sphere packing on $\Lambda$ of maximal density. 
    The \textbf{packing radius} $\rho$ of $\Lambda$ is the radius of the largest open ball in $\mathbb{R}^n$ whose translates by elements of $\Lambda$ do not overlap. 
\end{defin}
Let $E$ be an Elliptic curve over a field $k$ and let $E(k)$ denote the set of $k$-rational points of $E$. It is a well-known fact that there is a binary operation on $E(k)$ which makes it into an abelian group. Moreover, if $k$ is a number field or a function field over a finite field, then $E(k)$ is a finitely generated abelian group. Thus, we get that $E(k)_{\text{free}}:=E(k)/E(k)_{\text{torsion}}$ is a free abelian group of finite rank. There is also a positive-definite quadratic form $\hat{h}:E(k)_{\text{free}}\longrightarrow\mathbb{R}$ called the Néron–Tate height. The Néron–Tate height has an associated positive-definite symmetric bilinear form $\langle-,-\rangle:E(k)_{\text{free}}\times E(k)_{\text{free}}\longrightarrow \mathbb{R}$ given by $$\langle P,Q\rangle =\frac{1}{2}(\hat{h}(P+Q)-\hat{h}(P)-\hat{h}(Q)).$$
\begin{defin}
    A \textbf{lattice} is a finitely generated free abelian group equipped with a positive-definite symmetric bilinear form. The \textbf{Mordell-Weil} lattice of $k$-rational points of an elliptic curve $E$ is the group $E(k)_{\text{free}}$ equipped with the bilinear form described above. The sphere packing density associated to an elliptic curve $E$ is the density of the densest sphere packing on its Mordell-Weil lattice. 
\end{defin}
If $(\Lambda,\langle -,-\rangle )$ is a lattice where $\Lambda$ is of rank $n$, then we can find a group homomorphism $\Phi:\Lambda\longrightarrow \mathbb{R}^n$ such that $(\Phi(P))\cdot(\Phi(Q))=\langle P,Q\rangle $ for all $P$ and $Q$ in $\Lambda$. $\Phi$ can be constructed using the Gram-Schmidt process. Thus, instead of studying the lattice $(\Lambda,\langle -,-\rangle )$, we can study its image in $\mathbb{R}^n$ equipped with the usual dot product. Note that $\Phi$ is not unique. In fact, the images of $\Lambda$ under different $\Phi$ may be different. However, they will be isomorphic as lattices. Thus, we can define the $k$-density of the densest sphere packing associated to an elliptic curve $E$ as the density of the densest sphere packing of an embedding of the lattice $E(k)_\text{free}$ into $\mathbb{R}^n$.
\section{General Results}
\begin{prop}
Let $\rho$ be packing radius of a maximal rank abelian subgroup $\Lambda\subseteq \mathbb{R}^n$ and let $V(\Lambda)$ be the volume of the fundamental domain of $\Lambda$. Then the density of the densest sphere packing of $\Lambda$ is given by $$\frac{\pi^{\frac{n}{2}}}{V(\Lambda)\cdot\Gamma(\frac{n}{2}+1)}\rho^n.$$
\begin{proof}
    This follows from the formula for the volume of the $n$-ball of radius $\rho$ and the observation that the volume of the fundamental domain occupied by the spheres to achieve the densest packing is exactly the volume of an $n$-sphere of radius $\rho$. 
\end{proof}
\end{prop}
\begin{defin}
Let $\Lambda$ be a maximal rank free subgroup of $\mathbb{R}^n$. Then the \textbf{normalized center density} $\delta_\Lambda$ is given by 
$$
\frac{\rho^n}{V(\Lambda)}
$$ where $\rho$ is the packing radius of $\Lambda$ and $V(\Lambda)$ is the volume of the fundamental domain of $\Lambda$. 
\end{defin}
It is easy to see that for a fixed $n$, the ratio of the normalized sphere density and the sphere density is constant for all lattices in $\mathbb{R}^n$. Thus, having a large lower bound on the normalized sphere packing density gives us a large lower bound on the sphere packing density. 
\begin{prop} 
\cite{Elk}
Let $E$ be an elliptic curve over $K$ where $K$ is a number field or function field over a finite field. For any subfield $k$ of $K$, the normalized center density of $E(k)$ is given by $$
\Delta^{-1/2}\left(\frac{N_{\text{min}}}{4}\right)^{n/2}
$$ where $n$ is the rank of $E(k)$, $N_{min}$ is the minimal non-zero value of the Néron–Tate height function, and $\Delta$ is the absolute value of the determinant of the pairing matrix of the generators of $E(k)_{\text{free}}$. 
\begin{proof}
    Note that the volume of the fundamental domain of a lattice is the square root of the absolute value of the determinant of the Gram matrix of a set of generators of the lattice. To avoid overlaps in the Mordell-Weil lattice, the packing radius $\rho$ must be half the minimal distance between any two distinct points of the lattice. This minimal distance is given by $\sqrt{N_{min}}$.
\end{proof}
\end{prop}
\section{The Curve $E: y^2=x^3+t^q-t$}
Let $p>3$ be a prime such that $p\equiv-1$ mod $6$, $q$ an odd power of $p$, and $r$ a sufficiently large power of $p$ (the notion of ``sufficiently large" is made precise below). In this section, we consider the curve given by $y^2=x^3+t^q-t$ over the field $k:=\mathbb{F}_r(t)$. 
\begin{defin}
Fix an odd $p$-power $q=p^c$. We say that $r:=p^s$ is \textbf{sufficiently large} if $s$ is a multiple of $c$, $8$ divides $(p+1)s$, and $3(p^c-1)$ divides $p^s-1$. In particular, $\mathbb{F}_q\subseteq \mathbb{F}_r$. If $r$ is sufficiently large, then so are all powers of $r$.  
\end{defin}
\begin{defin}
The \textbf{naive height} $h$ of a $\mathbb{F}_r(t)$-rational point $P=(x,y)$ on an elliptic curve is $\text{deg}(x):=\text{max}(\text{deg}(f),\text{deg}(g))$ where $x=\frac{f(t)}{g(t)}$ for $f$ and $g$ coprime polynomials. If $P=[0:1:0]$, then we define $h(P)=0$.
\end{defin}
We wish to find a lower bound on $N_\text{min}$ for our elliptic curve. This is the same as finding the minimum value of $\hat{h}$ on non-torsion points. Since $\hat{h}$ has a bounded difference with the naive height, it is perhaps useful to find a lower bound on the naive height.
\begin{prop}
The naive height on $E$ is bounded below by $\frac{q+1}{3}$.
\begin{proof}
    Let $\text{ord}_\infty$ denote the valuation at $\infty$. Let $P=(x,y)$ be a point on $E$ that is not the identity. Then we can apply $\text{ord}_\infty$ to both sides of $y^2=x^3+t^q-t$ to get 
    $$
    2\text{ord}_\infty(y)=\text{ord}_\infty(x^3+t^q-t)
    .$$
Note that $\text{ord}_\infty(t)=-1\neq-q=\text{ord}_\infty(t^q)$. Thus, we may use the strict triangle inequality to get that $\text{ord}_\infty(t^q-t)=-q$. However, $q$ is not a multiple of $3$, whereas $\text{ord}_\infty(x^3)=3\text{ord}_\infty(x)$ is, so we must have, using the strict triangle inequality again, that $$\text{ord}_\infty(x^3+t^q-t)=\text{min}(3\text{ord}_\infty(x),-q).$$
We also know that $2\text{ord}_\infty(y)=\text{ord}_\infty(y^2)=\text{min}(3\text{ord}_\infty(x),-q)$. Since $2$ does not divide $q$, we must have that $3\text{ord}_\infty(x)< -q$. Thus, we get 
$$
\text{ord}_\infty(x)\leq -\left(\frac{q+1}{3}\right)
.$$
This shows that $x$ has a pole of order at least $\frac{q+1}{3}$ at $\infty$, so we must have $h(P)=\text{deg}(x)\geq \frac{q+1}{3}$.
\end{proof}
\end{prop}
Using this proposition, we can now prove that the naive height and the Néron–Tate height are the same for our curve.
\begin{prop}
The Néron–Tate height and the naive height agree on $E$.
\begin{proof}
    Since the Neron-Tate height $\hat{h}$ is uniquely characterized by the fact that its difference with the naive height $h$ is bounded and $\hat{h}(2P)=4\hat{P}$ for all $\mathbb{F}_r(t)$-rational points $P$ on the curve, it suffices to show that the naive height on $E$ satisfies 
    $$
    h(2P)=4h(P)
    $$
    for all rational points points $P$ on the curve $E$. Clearly, this holds for when $P$ is the identity. Let $P=(x,y)$ be a non-identity $\mathbb{F}_r(t)$-rational point on the curve where $x=\frac{f(t)}{g(t)}$ where $f$ and $g$ are coprime polynomials. By the proof of Proposition 4.3, we have
    $$\text{deg}(g)-\text{deg}(f)=\text{ord}_\infty(x)\leq -\left(\frac{q+1}{3}\right).$$
    This implies that $\text{deg}(f)\geq \text{deg}(g)+\frac{q+1}{3}$, so we conclude that $\text{deg}(x)=\text{deg}(f)$ for all points $P=(x,y)$ with $x=\frac{f}{g}$ where $f$ and $g$ are coprime polynomials.\\
    Now, we compute $2P$ explicitly. Using basic arithmetic of elliptic curves, we see that the first coordinate of the point $2P$ is 
    $$
    \frac{f^4-8fg^3\cdot(t^q-t)}{4(gf^3+g^4\cdot(t^q-t))}
    .$$
    By our earlier observation on the height of $\mathbb{F}_r(t)$-rational points of $E$, if there is no cancellation between the numerator and the denominator, the naive height of $2P$ is $\text{deg}(f^4-8fg^3\cdot(t^q-t))$. However, using Proposition 4.3, we get  $$\text{deg}(f^4)=4\text{deg}(f)\geq \text{deg}(f)+3\text{deg}(g)+q+1>\text{deg}(8fg^3\cdot(t^q-t)).$$
    Thus, we get $\text{deg}(f^4-8fg^3\cdot(t^q-t))=\text{deg}(f^4)=4\text{deg}(f)=4h(P)$.\\
    All we need to show now is that 
    $$
    \frac{f^4-8fg^3\cdot(t^q-t)}{4(gf^3+g^4\cdot(t^q-t))}
    $$
    has no cancellation. Assume that $\tau$ is an irreducible (in $\mathbb{F}_r[t]$) that divides both the numerator and denominator. Then $\tau$ must also divide
    $$
    g(f^4-8fg^3\cdot(t^q-t))-\frac{f}{4}(4(gf^3+g^4\cdot(t^q-t)))=-9g^4f\cdot(t^q-t).
    $$
    Since $\tau$ is irreducible and $9\neq 0$, we must then have that $\tau$ either divides $g,f,$ or $t^q-t$. \\
    \\
    \textbf{Case 1:} $\tau$ divides $g$\\
    Since $\tau$ divides both $f^4-8fg^3\cdot(t^q-t)$ and $g$, we must have that it divides $f$ as well. This is a contradiction as $f$ and $g$ were assumed to be coprime polynomials.\\
    \textbf{Case 2:} $\tau$ divides $f$\\
    Since $\tau$ divides both $gf^3+g^4\cdot(t^q-t)$ and $f$, we must have that it divides either $g$ or $t^q-t$. Since $g$ and $f$ are coprime polynomials, this can only happen if $\tau$ divides $t^q-t$. However, $\mathbb{F}_q\subseteq \mathbb{F}_r$ and the roots of $t^q-t$ are precisely the elements of $\mathbb{F}_q$. Thus, the irreducible $\tau$ must be of the form $t-\alpha$ for some $\alpha\in\mathbb{F}_q$. This implies that $\alpha$ is a root of $f$. Since $x=\frac{f(t)}{g(t)}$ where $f$ and $g$ are coprime polynomials, we also get that $x$ has a root at $\alpha$. Note that $t^q-t$ has a root at $\alpha$ of multiplicity one, so we can use the strict triangle inequality to get 
    $$
    2\text{ord}_\alpha(y)=\text{ord}_\alpha(y^2)=\text{ord}_\alpha(x^3+t^q-t)=\text{min}(\text{ord}_\alpha(x^3),\text{ord}_\alpha(t^q-t))=1,
    $$
    which is a contradiction as $2$ does not divide $1$.
    \\
    \textbf{Case 3:} $\tau$ divides $t^q-t$\\
    Since $\tau$ divides both $f^4-8fg^3\cdot(t^q-t)$ and $t^q-t$, we must have that it divides $f$ as well. This then leads to the same contradiction as in Case 2.
    \\
    This shows that there is no cancellation in $\frac{f^4-8fg^3\cdot(t^q-t)}{4(gf^3+g^4\cdot(t^q-t))}$, so the naive height satisfies $h(2P)=4h(P)$.
\end{proof}
\end{prop}
\begin{rem}Using Proposition 4.4, we can give an elementary argument to show that $E(\mathbb{F}_r(t))$ is torsion-free for $r$ sufficiently large.
Let $P=(x,y)$ be a non-trivial $\mathbb{F}_r(t)$-rational point on $E$ that is torsion. By Proposition 4.4, we know that $0=\hat{h}(P)=\text{deg}(x)$. This is only possible if $x=c$ for some $c\in\mathbb{F}_r$. This would imply that $y^2=t^q-t+c^3$. This means that $y$ is a polynomial in $t$. Since $q$ is an odd integer, $t^q-t+c^3$ cannot be the square of a polynomial. This is a contradiction. 
\end{rem}
\begin{prop}
[Proposition 8.4.1(3) in \cite{GU}] For $r$ sufficiently large, rank$(E(\mathbb{F}_r(t)))=2(q-1)$.
\end{prop}
\begin{prop}
[Corollary 9.2(3) in \cite{GU}] $\Delta$ is bounded above by $r^{\lfloor \frac{q}{6} \rfloor}$ for $r$ sufficiently large.
\end{prop}
\begin{prop}
For $r$ sufficiently large, the normalized center density $\delta_{E(\mathbb{F}_r(t))}$ is bounded below by
$$\frac{\sqrt{|\Sha(E/\mathbb{F}_r(t))|}}{(r^{\lfloor \frac{q}{6} \rfloor})^{1/2}}\cdot\left(\frac{q+1}{12}\right)^{q-1}.$$
\begin{proof}
This follows by combining the results of Proposition 3.3, Proposition 4.3, Proposition 4.4, Proposition 4.5, and Proposition 4.6. 
\end{proof}
\end{prop}
Using the trivial bound on the Tate–Shafarevich group, we can now use Proposition 4.7 to get explicit lower bounds on the densities of our family of curves. However, as we will see, this lower bound can often be improved by considering sublattices of the Mordell-Weil lattice. Indeed, let $\Lambda$ be a maximal rank sublattice of the Mordell-Weil lattice of $E$. We may then replace $\Delta$ by $V(\Lambda)$ in the in Proposition 3.3 to find a lower bound on the normalized packing density of E. With this as the motivation, we list some explicit rational points on the curve and consider the sublattice that they generate. 
\begin{prop}
Let $\sigma$ be a solution to $\sigma^{6(q-1)}=-1$ and $\beta$ a solution to $\beta^q+\beta=1$ in $\mathbb{F}_r$ for $r$ sufficiently large. Then $P=(\sigma^2(t-(\beta/\sigma^6))^\frac{q+1}{3},\sigma^3(t-(\beta/\sigma^6)^q)^\frac{q+1}{2})$ is a point in $E(K)$. Thus, we get $6q(q-1)$ $\mathbb{F}_r(t)$-rational points on $E$.
\begin{proof}
    This can be verified by plugging in the point $P$ into the defining equation of our elliptic curve and using the relations on $\sigma$ and $\beta$. Alternatively, these points can be deduced by considering the map $C_{6,q}\longrightarrow E_0$ which we get by presenting $C_{6,q}$ and $E_0:y^2=x^3+1$ as quotients of the Fermat curve of degree $q+1$. 
\end{proof}
\end{prop}
Using \MYhref{https://github.com/Arjun-Nigam/MAGMA-code-lattices/blob/main/MAGMA\%20code\%20for\%20t\%5Eq\%20-\%20t\%20curve.txt}{this MAGMA script}, we find the volume of (an Euclidean embedding of) the sublattice generated by the points in Proposition 4.8 which we then use to find lower bounds on the sphere packing density of the Mordell-Weil lattice for different values of $p,q,$ and $r$. We compare it with the lower bounds obtained using Proposition 4.7 and the trivial bound $|\Sha(E/\mathbb{F}_r(t))|\geq 1$ in the table below.
\begin{center}
\begin{tabular}{ | m{1.2em} | m{1.2em}| m{2.1cm} | m{3.5cm} | m{1.5cm} | m{2.5cm} | } 
  \hline
  $q$ & $r$ & Lower Bound on Normalized Density (Using MAGMA) & Lower Bound on Normalized Density (Using Proposition 4.7 and the trivial bound $|\Sha(E/\mathbb{F}_r(t))|\geq 1$) & Dimension & Best Known (Normalized) Sphere Packing Density (If Known) \cite{Coh}\\ 
  \hline
  $5$ & $5^4$ & $0.0625$ & $0.0625$ & $8$ & $0.0625$\\ 
  \hline
  $5$ & $5^8$ & $0.0625$ & $0.0625$ & $8$ & $0.0625$\\ 
  \hline
  $5$ & $5^{12}$ & $0.0625$ & $0.0625$ & $8$ & $0.0625$\\ 
  \hline
  $5^3$ & $5^{16}$ & $\sim 2.653\times10^{20}$ & $\sim 6.198\times 10^{12}$ & $248$ & \\ 
  \hline
  $11$ & $11^2$ & $\sim 0.0909$ & $\sim 0.0909$ & $20$ & $\sim 0.1315$\\ 
  \hline
  $11$ & $11^6$ & $\sim 0.0909$ & $\sim 0.00075$ & $20$ & $\sim 0.1315$\\ 
  \hline
  $17$ & $17^4$ & $\sim 2.272$ & $\sim 0.0078$ & $32$ & $\sim 2.565$\\ 
  \hline
\end{tabular}
\end{center}

We see that there is repetition in the lower bound of the normalized density (computed using MAGMA) when $q$ is fixed and $r$ varies. This is to be expected as the explicit points in Proposition 4.8 and the explicit formulation of the Néron–Tate height do not depend on $r$ (assuming that $r$ is sufficiently large). Thus, they generate the same lattice. The same repetition is not observed in the lower bound obtained using Proposition 4.7. This is likely due to the fact that we used the trivial bound on $|\Sha(E(\mathbb{F}_r(t)))|$ when computing it and $\Sha(E(\mathbb{F}_r(t)))$ may depend on $r$.

Up to this point, we have been using properties of the given elliptic curves to deduce information about the densities of their Mordell-Weil lattices. Now, we can use the two lower bounds to say something interesting about the elliptic curve $E$. The lower bound obtained using MAGMA divided by the lower bound obtained using Proposition 4.7 is a lower bound on $\sqrt{\Sha(E/\mathbb{F}_r(t))}$. This follows from how we computed these lower bounds. In particular, we see that for the cases considered above, $E(E/\mathbb{F}_r(t))$ is non-trivial when $q=5^3$ and $r=5^{16}$ or when $q=17$ and $r=17^4$.

\begin{rem}$E_8$ is the unique lattice (up to isometry and rescaling) to have the highest density sphere packing in dimension $8$. For $q=5$ and $r=5^4$, we observe that the lower bound obtained equals the sphere packing density of $E_8$. This means that the Mordell-Weil lattice must have sphere packing density equal to $E_8$ which, by uniqueness, implies that the Mordell-Weil lattice is $E_8$. This also implies that $E(E/\mathbb{F}_r(t))$ is trivial in this case.
\end{rem}
\section{The Legendre Curve}
Let $p$ be an odd prime and $d=p^f+1$ for some positive integer $f$. Consider the function field $K_d:=\mathbb{F}_p(\mu_d,u)$ where $\mu_d$ is the set of primitive $d^\text{th}$ roots of unity and $u^d-t=0$ for indeterminate $t$. We study the elliptic curve $$
E:y^2=x(x+1)(x+t).
$$
As with the previous curve, we would like to have explicit points on this curve.
\begin{prop}
Let $E$ and $K_d$ be as defined. Then for a fixed primitive $d^{\text{th}}$ root of unity $\zeta$ and $i\in\{0,1,\ldots,d-1\}$, the point given by $
P_i^{(d)}:=(\zeta^i_du,\zeta^i_du(\zeta^i_du+1)^{d/2})
$
is a $K_d-$rational point of $E$.
\begin{proof}
    This can be checked by plugging in each $P^{(d)}_i$ into the equation defining $E$ or by observing the action of $Gal(K_d/K)$ on the $K_d-$ rational point $P:=(u,u(u+1)^{d/2})$.
\end{proof}
\end{prop}
\begin{prop}
The set of explicit rational points in Proposition 5.1 generate a rank $d-2$ subgroup of the Mordell-Weil group. $E(K_d)$ also has rank $d-2$.
\begin{proof}
The first part of this proposition follows from Corollary 4.3 in \cite{Ulm} and the second part follows from Corollary 5.3 of the same paper.
\end{proof}
\end{prop}
Thus, we get that the subgroup generated by the $P_i^{(d)}$ is of finite index in $E(K_d)$.
Just like in the previous section, it is difficult to find a set of generators for the free part of the Mordell-Weil group. However, instead of finding a lower bound on the normalized sphere packing density of $E$, we can find a lower bound on the normalized sphere packing density of the finite index sublattice generated by the $P_i^{(d)}$. 
\begin{prop}
The height pairing of the points $P^{(d)}_i$ is given by
$$
  \langle P^{(d)}_i,P^{(d)}_j\rangle  = \begin{dcases}
        \frac{(d-1)(d-2)}{2d} & \text{if } i=j \\
        \frac{1-d}{d} & \text{if } i\neq j \text{ and } i-j \text{ is even} \\
       0 & \text{if } i-j \text{ is odd} \\
    \end{dcases}
$$
\begin{proof}
This is Theorem 8.2 in \cite{Ulm}
\end{proof}
\end{prop}
\begin{prop}
For any non-torsion point $Q$ in the subgroup generated by the $P_i^{(d)}$, we have the inequality $$\frac{d-1}{2d}\leq \hat{h}(Q).$$
\begin{proof}
Let $Q=\sum_{0}^{d-1} a_iP_i^{(d)}$ be a point in the subgroup generated by the $P_i^{(d)}$. Then we can use bilinearity of the height pairing and the previous proposition to get $$\hat{h}(Q)=\left\langle \sum_{0}^{d-1} a_iP_i^{(d)},\sum_{0}^{d-1} a_iP_i^{(d)}\right\rangle =\frac{(d-1)(d-2)(a_0^2+\cdots+a_{d-1}^2)+(1-d)(\sum_{i-j\text{ is even}}2a_ia_j)}{2d}.$$
Thus, we may write $\hat{h}(Q)$ as $\frac{d-1}{2d}\cdot m$ where $m$ is some integer. Since the height function is positive definite, we must have that $m\geq 1$. This gives us the required lower bound. 
\end{proof}
\end{prop}
Using these results, we can now compute explicit lower bounds on the normalized sphere packing density of a sublattice of $E$ using \MYhref{https://github.com/Arjun-Nigam/MAGMA-code-lattices/blob/main/MAGMA\%20code\%20for\%20Legendre\%20curve.txt}{this MAGMA script}.\\
\begin{center}
\begin{tabular}{ | m{2em} | m{1cm}| m{1.5cm} | m{3cm} | m{4cm} |} 
  \hline
  $p$ & $f$ & Dimension & Normalized Density (Lower Bound) & Best Known (Normalized) Sphere Packing Density \cite{Coh}\\
  \hline
  $3$ & $1$ & $2$ & $\sim 0.125$ & $\sim 0.288$ \\ 
  \hline
  $3$ & $2$ & $8$ & $\sim 1.953\times 10^{-6}$ & $\sim 0.0625$ \\ 
  \hline
  $3$ & $3$ & $26$ & $\sim 3.208\times 10^{-26}$ & $\sim 0.577$\\
  \hline
  $5$ & $1$ & $4$ & $\sim 0.005$ & $\sim 0.125$\\ 
  \hline
  $5$ & $2$ & $24$ & $\sim 8.119 \times 10^{-24}$ & $\sim 1.003$\\ 
  \hline
  $7$ & $1$ & $6$ & $\sim 1.22\times 10^{-4}$ & $\sim 0.0721$\\
  \hline
\end{tabular}
\end{center}
We observe that as $p$ and $f$ get larger, the normalized density gets smaller. This is because the Gram matrix obtained from the height pairing has a large determinant for large values of $p$ and $f$.

Note that, unlike the previous section, the computational lower bounds that we obtained here are not necessarily lower bounds on the density of $E$. While we were able to obtain a lower bound on the minimum non-zero value of the Néron–Tate height function on the sublattice generated by the $P_i^{(d)}$, we do not know if this value is also the minimum on the entire Mordell-Weil lattice. 
\section{Generalizations}
We now discuss how the methods used in this paper may be generalized to other ``nice" families of elliptic curves. Our strategy in Section 4 is motivated by \cite{Elk}. We first proved that the Néron–Tate height and the naive height agree for our family of elliptic curves and then found a lower bound for the naive height using elementary techniques. While the proofs of these results use properties specific to the family of curves being studied, there are many examples in the literature of elliptic curves where the naive and Néron–Tate height agree (for example, the family of hyperelliptic curves considered in \cite{Elk}). 

Even if the naive and Néron–Tate heights do not agree for all points, we might still be able to perform explicit computations on the lower bounds of lattice densities. Let $E$ be an elliptic curve over the field $\mathbb{F}_q(t)$ where $q$ is a $p$-power for a prime $p>3$. We choose a Weierstrass model 
$$
y^2 = x^3 + a_4x + a_6
$$
for $E$ where $a_i\in \mathbb{F}_q[t]$ are polynomials and max$(\lceil \frac{\text{deg}(a_i)}{i}\rceil)$ is minimal (doing so is possible by lemma 5.43 of \cite{Shi}). We define $d$ to be the smallest integer such that $a_i\leq di$. We say that $E$ is a ``nice" elliptic curve if 
\begin{enumerate}
    \item For all places $\nu$, the discriminant $\Delta$ has no roots, simple roots, or double roots at $\nu$.
    \item If $\Delta$ has double roots at $\nu$, then the reduction of $E$ at $\nu$ is additive.
\end{enumerate}
Assume that $E$ is a ``nice" curve and let $P$ be a $\mathbb{F}_q(t)$-rational point. We have $\langle P,P\rangle=2\chi + 2(P\cdot O)+\sum_v \text{contr}_v(P)$ where the sum is taken over all singular reducible fibers and $P\cdot O$ is the intersection number (Theorem 6.24 in \cite{Shi}). Since we assumed that $E$ is ``nice", table 15.1 of \cite{Silver} tells us that all singular fibers are irreducible. Thus, we get $\langle P,P\rangle=2\chi + 2(P\cdot O)$. However, since $\chi=d$, all we need to do to compute the height is to find the intersection number $P\cdot O$. 

Let $(x_0(t),y_0(t))$ be the coordinate representation of the point $P$ where $x_0=\frac{f(t)}{g(t)}$ for polynomials $f$ and $g$. Since $a_i$ are polynomials, they do not blow up when $t$ is finite. Thus, given a finite place $t_0$, we have that $x_0(t_0)$ blows up if and only if $y_0(t_0)$ blows up. This means us that the section determined by $P$ meets the section determined by $O$ at deg$(g)$-many times (counted with multiplicity). 

For $t=\infty$, we change coordinates of our curve by $x=t^{2d}x'$ and $y=t^{3d}y'$. By our choice of $a_i$, we get that the coefficients in our new model are regular functions at $t=\infty$. Thus, we get that the section determined by $P$ intersects the section determined by $O$ at infinity max$(0,\text{deg}(f)-\text{deg}(g)-2d)$-many times. If deg$(f)$ large enough, we get that the section determined by $P$ meets the section determined by $O$ deg$(g)+\text{deg}(f)-\text{deg}(g)-2d)=\text{deg}(f)-2d$ many times. Accounting for the degree $2$ morphism $E\longrightarrow \mathbb{P}^1$, we get $P\cdot O=\frac{1}{2}\left(\text{deg}(f)-2d\right)$. This gives us $\langle P,P\rangle=2d + 2(P\cdot O)=\text{deg}(f)$. 

We have now shown that for ``nice" curves $E$ and points $P=(\frac{f}{g},y)$ with $\text{deg}(f)$ is large enough, the naive and Néron–Tate heights on $P$ agree. Thus, if we have a lower bound on the naive height for $E$ and a set of points $P_i=(\frac{f_i}{g_i},y_i)$ where deg$(f_i)$ is large enough for each $i$ and the $P_i$'s generate a maximal rank sublattice $L$, we can compute lower bounds on the lattice density of $L$. 
\section{Acknowledgements}
I would like to thank Dr. Douglas Ulmer for suggesting the project and mentoring me with it. Dr. Ulmer's constant encouragement, expertise, and patient guidance were invaluable in the writing of this paper. I would also like to thank Dr. Bryden Cais, Dr. Steven J. Miller, Daniel Lewis, and Tristan Phillips for their valuable feedback on the first draft of this paper.
\printbibliography
\end{document}